\documentclass{birkjour}
\usepackage{mathtools,commath}
\usepackage{amssymb,amsmath,amsthm}
\usepackage{enumitem}
\usepackage{bbm}
\usepackage{hyperref}

\usepackage{bigints}

\newtheorem{theorem}{Theorem}
\newtheorem{lemma}[theorem]{Lemma}

\newtheorem{remark}[theorem]{Remark}
\newtheorem{cor}[theorem]{Corollary}

\newtheorem{example}[theorem]{Example}
\numberwithin{equation}{section}
\numberwithin{theorem}{section}
\renewenvironment{proof}{\hspace*{-0.6cm}{\bfseries Proof:}}{\hfill $\blacksquare$}
\def\bR{\mathbb R}

\begin{document}

\title[Nonexistence of positive entire solutions]{Nonexistence of positive entire solutions to systems of semilinear elliptic inequalities }

\author{ Daniel Devine} 
\address{School of Mathematical Sciences, University College Cork, Cork, Ireland}
\email{ ddevine@ucc.ie} 

\begin{abstract}
Consider the system
\begin{equation*}
\left\{
\begin{aligned}
&\Delta u \geq p(x)g(v)   \quad&&\mbox{ in } \bR^n, \\
&\Delta v \geq q(x)f(|\nabla u|) \quad&&\mbox{ in } \bR^n,\\
&u>0, v>0 \quad&&\mbox{ in } \bR^n,
\end{aligned}
\right.
\end{equation*}
where $n\geq2$, $f,g\in C[0,\infty)$ and $p,q\in C(\bR^n)$. If $f,g$ are non-decreasing and convex, using a comparison argument we prove a sufficient Keller-Osserman type condition for the nonexistence of entire solutions. This result is sharp for a wide class of systems. Similar results are also shown to hold for systems without gradient terms.
\end{abstract}
\maketitle
\begingroup
\renewcommand{\thefootnote}{}
\footnotetext{This version of the article has been accepted
for publication, after peer review (when applicable) but
is not the Version of Record and does not reflect
post-acceptance improvements, or any corrections.
The Version of Record is available online at:
\url{https://doi.org/10.1007/s00030-026-01258-4}.}
\endgroup
\textbf{Mathematics Subject Classification (2020):} 35B44, 35J47, 35J91

\textbf{Keywords:} elliptic systems, coercive systems, nonlinear gradient terms

\section{Introduction}
In this paper, we study the nonexistence of entire solutions to the following coercive semilinear elliptic system of inequalities
\begin{equation}\label{system}
\left\{
\begin{aligned}
&\Delta u \geq p(x)g(v)   \quad&&\mbox{ in } \bR^n, \\
&\Delta v \geq q(x)f(|\nabla u|) \quad&&\mbox{ in } \bR^n,\\
&u>0, v>0 \quad&&\mbox{ in } \bR^n,
\end{aligned}
\right.
\end{equation}
where $n\geq2$,  $f,g:[0,\infty)\to [0,\infty)$ are continuous and positive on $(0,\infty)$, while $p,q\in C(\bR^n)$ are positive on $\bR^n\setminus\{0\}$.  By a solution of \eqref{system}, we shall always mean a pair $(u,v)\in C^2(\bR^n)^2$ which satisfies \eqref{system} pointwise. The goal of this paper is to provide what appear to be the first results dealing with non-radial solutions of systems of the form \eqref{system}. Semilinear coercive systems with gradient nonlinearities such as \eqref{system} have received considerable attention in the last two decades \cite{BGW22,D2025,DK24,DS23,FV17,GGS19, S15}, but all the work to date has been focused on solutions which are radially symmetric. We note that for scalar problems with gradient nonlinearities there is an incredibly rich literature, and we point the reader to the classical papers \cite{AGQ2012, FS2011, FPR2010} and the references therein. 

The prototype system of interest is
\begin{equation}\label{DLSsystem}
\left\{
\begin{aligned}
&\Delta u = v &&\mbox{in}\ \Omega, \\
&\Delta v = |\nabla u|^2&&\mbox{in}\ \Omega,
\end{aligned}
\right.
\end{equation}
which was introduced by D\'{i}az \textit{et al.} \cite{DLS2005}. This corresponds to system \eqref{system} for $p\equiv q\equiv 1$,  $g(t)=t$ and $f(t)=t^2$, in the case of equality. System \eqref{DLSsystem} arose from studying steady states of the corresponding parabolic system, 
\begin{equation*}
\left\{
\begin{aligned}
&u_t-\Delta u = v &&\mbox{in}\ \Omega\times(0,T), \\
&v_t-\Delta v = |\nabla u|^2 &&\mbox{in}\ \Omega\times(0,T),
\end{aligned}
\right.
\end{equation*}
which was used to model the dynamics of a fluid in situations where viscous heating effects cannot be neglected. Considering a unidirectional flow, independent of distance in the flow direction, the speed $u$ and temperature $v$ satisfy the above coupled equations, where the source terms $v$ and $|\nabla u|^2$ represent the  buoyancy force and viscous heating, respectively. It is conjectured that solutions of this parabolic system will experience finite time blow-up in a suitable norm, analogous to classical Fujita type results for reaction-diffusion equations and systems. There is a huge literature on such problems, so we point the reader to the recent article \cite{FilLom20} and the survey article \cite{Soup01} for results concerning parabolic equations with gradient nonlinearities, and to \cite{EH91, QS07} for results on parabolic systems.

The authors of \cite{DLS2005} studied the system \eqref{DLSsystem}  in an open ball $B_R$ centered at the origin. The main result of \cite{DLS2005} states that for each $R>0$, there exists a unique radial solution $(u,v)$ of \eqref{DLSsystem} satisfying $u(0)=0, v(0)>0$ and $|(u(x),v(x))|\to\infty$ as $|x|\to R^{-}$. Moreover, the map $v(0)\mapsto R$ is continuous and strictly decreasing. Note that since system \eqref{system} depends only on $\nabla u$, for any solution $(u,v)$ and $c\in\bR$, $(u+c,v)$ is also a solution.

The results of \cite{DLS2005} were later extended by Singh \cite{S15}, who derived optimal conditions for the existence of radial solutions of \eqref{system} (in the case of equations instead of inequalities) for $p\equiv q\equiv 1$, $g(t)=t^k$ for $k>0$. In particular, he showed that for the case $p\equiv q\equiv 1$, $g(t)=t^k$, system \eqref{system} admits radially symmetric solutions if and only if
\begin{equation}\label{KOcond}
      \int_{1}^{\infty}\frac{ds }{\left(\int_{0}^{s}\sqrt{f(t)}dt\right)^{2k/(2k+1)}}\ ds=\infty,
\end{equation}
a condition which is clearly not satisfied for system \eqref{DLSsystem}.
Singh's work therefore elucidated the fact that system \eqref{DLSsystem} does not admit any entire radial solutions, but left open the question of non-radial solutions. We note that a condition such as \eqref{KOcond} is known as a Keller-Osserman type condition, so named for the seminal work independently performed by Keller \cite{K57} and Osserman \cite{O57} on the problem
\begin{equation*}
\left\{
\begin{aligned}
\Delta u&=f(u)  &&\quad\mbox{ in } \Omega, \\
u&\to\infty &&\quad\mbox{ as } x\to\partial\Omega,
\end{aligned}
\right.
\end{equation*}
where $f\in C^{1}([0,\infty))$ and $\Omega\subset\bR^n$ is bounded and smooth. Following Singh's work \cite{S15}, several variations and generalisations of system \eqref{system} have been investigated, see for instance \cite{ BGW22,D2025,DK24,DS23,FV17,GGS19}. We emphasise again that in all cases, these studies focus solely on radially symmetric solutions, and do not address non-radial solutions. To the best of our knowledge the results of this paper are the first to do so. Moving forward, we shall assume one or more of the following conditions: \\
\begin{itemize}
    \item[(A1)] $p,q\in C(\bR^n)$ are positive on $\bR^n\setminus\{0\}$. $f,g$ are continuous on $[0,\infty)$ and positive on $(0,\infty)$. $f$ is non-decreasing and $g$ is strictly increasing,
\\
    \item [(A2)] there exist radially symmetric functions $\Tilde{p}=\Tilde{p}(|x|)$ and $\Tilde{q}=\Tilde{q}(|x|)$ satisfying
    \begin{equation*}
        p(x)\geq \Tilde{p}(|x|)\ \ \mbox{and}\ \ q(x)\geq \Tilde{q}(|x|),
    \end{equation*}
    where $\tilde{p},\tilde{q}$ are continuous on $[0,\infty)$ and positive on $(0,\infty)$.
\\    
    \item [(A3)] $f,g$ are convex.
\end{itemize}
\medskip

Our first main result is the following.
\begin{theorem}\label{main thm}
Assume (A1)-(A3), and suppose \eqref{system} admits a positive entire solution. Then the following system \normalfont{
\begin{equation}\label{systemrad}
\left\{
\begin{aligned}
&\Delta u = \Tilde{p}(|x|)g(v)   \quad&&\mbox{ in } \bR^n, \\
&\Delta v = \Tilde{q}(|x|)f(|\nabla u|) \quad&&\mbox{ in } \bR^n,\\
&u>0, v>0 \quad&&\mbox{ in } \bR^n,
\end{aligned}
\right.
\end{equation}
}\textit{admits a radially symmetric positive entire solution}.
\end{theorem}
Theorem \ref{main thm} says that, in order to prove the nonexistence of solutions to \eqref{system} (radial or non-radial), it is enough to know whether or not radial solutions of \eqref{systemrad} exist. At this point we note that it is in fact possible to treat systems where the monotonicity and convexity assumptions (A1) and (A3) do not hold, although we still need to compare with a system where they do. Explicitly, if instead of system \eqref{system} we consider
\begin{equation}\label{sysgeneral}
\left\{
\begin{aligned}
&\Delta u \geq p(x)\phi(v)   \quad&&\mbox{ in } \bR^n, \\
&\Delta v \geq q(x)\psi(|\nabla u|) \quad&&\mbox{ in } \bR^n,\\
&u>0, v>0 \quad&&\mbox{ in } \bR^n,
\end{aligned}
\right.
\end{equation}
where $\phi,\psi\in C[0,\infty)$ satisfy $\phi(t)\geq g(t)$ and $\psi(t)\geq f(t)$ for all $t\geq 0$, with $f,g,p,q$ satisfying (A1)-(A3), then a positive entire solution to system \eqref{sysgeneral} is automatically a positive entire solution to system \eqref{system}. The existence of a radially symmetric positive entire solution to \eqref{systemrad} is thus a necessary condition for the existence of a positive entire solution to system \eqref{sysgeneral}. For simplicity, most of our results are stated directly for system \eqref{system}.
\begin{remark}
\normalfont{
Results analogous to Theorem \ref{main thm} are well known in the scalar case. Indeed, for scalar problems far more general quasilinear elliptic inequalities such as
\begin{equation*}
    \mbox{div}(A(|\nabla u|)\nabla u)\geq f(x,u,\nabla u)\quad\mbox{in}\ \bR^n
\end{equation*}
may be treated. For example, in \cite{NU1997} Naito and Usami derive nonexistence results for the problem
\begin{equation*}
     \mbox{div}(A(|\nabla u|)\nabla u)\geq f(u)\quad\mbox{in}\ \bR^n.
\end{equation*}
By considering the corresponding radial equation in a manner analogous to Theorem \ref{main thm}, the authors obtain nonexistence results under quite general assumptions on the differential operator and nonlinear term $f$.
}
\end{remark}
We now move on to analysing the radially symmetric solutions of system \eqref{systemrad} in more detail. Before stating our results, it will be convenient to define
\begin{equation}\label{frakG}
    \mathfrak{G}(s)=\int_{0}^{s}g^2(t)\ dt\quad\mbox{for}\ s\geq 0,
\end{equation}
and we also define $F(s)=\int_{0}^{s}f(t)dt$ and $\mathcal{F}(s)=\int_{0}^{s}F(t)dt$. In order to directly study the nonexistence of radial solutions of system \eqref{systemrad}, we will need to impose some further conditions on the weight functions $\tilde{p},\tilde{q}$. While these technical conditions will be needed in the direct study of system \eqref{systemrad}, through the use of comparison arguments we can treat a very general class of weight functions, as we discuss below. We introduce the notation
\begin{equation}\label{phi}
    \Phi_{h}(r)\coloneqq \Phi (h,r)=h(r)-\frac{n-1}{r^n}\int_{0}^{r}t^{n-1}h(t)\ dt,\quad r>0,
\end{equation}
where $h(r)$ is any $C[0,\infty)$ function. We require that $\tilde{p}$ satisfies
\begin{equation} \label{inf p}
   \Phi_{\tilde{p}}(r)>0\quad\mbox{for all}\ r>0\ \mbox{and}\ \int^{\infty}\Phi_{\tilde{p}}(r)\ dr=\infty,
\end{equation}
where the lower bound of integration is omitted to allow any positive number, and we require that $\tilde{q}$ satisfies
\begin{equation}\label{inf q}
\mbox{for all}\ \rho>0\ \mbox{there exists}\ C_{\rho}>0\ \mbox{with}\ \inf\limits_{r\geq \rho}\Phi_{\tilde{q}}(r)\geq C_{\rho}.
\end{equation}
We note that any function satisfying condition \eqref{inf q} also satisfies condition \eqref{inf p}. Finally, we will assume that
\begin{equation}\label{p over q bound}
   \tilde{p}(r)\ \mbox{is bounded on}\ [0,\infty).
\end{equation}
 The prototype class of weight functions $\tilde{p},\tilde{q}$ considered for this analysis were functions which are non-decreasing on $[0,\infty)$. Conditions \eqref{inf p} and \eqref{inf q} are easily seen to hold if this is the case. Besides non-decreasing functions, a simple calculation shows that for any radially symmetric function $h\in C(\bR^n)$ satisfying
 \begin{equation*}
     \frac{d}{dr}\left(\frac{r}{|B_r|}\int_{B_r}h(x)\ dx\right)> 0\quad\mbox{for all}\ r>0,
 \end{equation*}
 one has $\Phi_h(r)>0$ for all $r>0$, using the notation $h(x)=h(|x|)$. Here, $|B_r|$ denotes the volume of the open ball of radius $r>0$ centred at the origin. Our first result concerning the positive radial solutions of system \eqref{systemrad} is the following.
\begin{theorem}\label{main thm 2}
Assume (A1) and (A2). Suppose further that $\Tilde{p}, \Tilde{q}$ satisfy \eqref{inf p}-\eqref{p over q bound}. If for every positive constant $C>0$ the Keller-Osserman type condition
\begin{equation}\label{KOgen}
    \int_{}^{\infty}\frac{ds}{g\left(\mathfrak{G}^{-1}\left(C\mathcal{F}(s) \right)\right)} <\infty
\end{equation}
holds, then system \eqref{systemrad} does not admit a radially symmetric positive entire solution. Consequently, if (A3) is satisfied and \eqref{KOgen} holds, then system \eqref{system} does not admit any positive entire solutions.
\end{theorem}
In condition \eqref{KOgen} we have once again omitted the lower bound of integration to allow any positive number.  While conditions \eqref{inf p}-\eqref{p over q bound} are needed to directly analyse system \eqref{systemrad} and derive condition \eqref{KOgen}, via a comparison argument we can show the same conclusion holds for systems with far more general weight functions, and more general nonlinear terms as discussed previously. Our next result demonstrates this.
\begin{cor}\label{corollary}
Assume that $f,g,\tilde{p},\tilde{q}$ satisfy the conditions of Theorem \ref{main thm 2}. Suppose that $\alpha=\alpha(|x|)$ and $\beta=\beta(|x|)$ are $C[0,\infty)$ functions satisfying $\alpha\geq \tilde{p}$ and $\beta\geq \tilde{q}$ at all points, and that $\phi,\psi$ are $C[0,\infty)$ functions satisfying $\phi\geq g$ and $\psi\geq f$ at all points. Consider the system
\begin{equation}\label{systemrad1}
\left\{
\begin{aligned}
&\Delta u = \alpha(|x|)\phi(v)   \quad&&\mbox{ in } \bR^n, \\
&\Delta v = \beta(|x|)\psi(|\nabla u|) \quad&&\mbox{ in } \bR^n,\\
&u>0, v>0 \quad&&\mbox{ in } \bR^n.
\end{aligned}
\right.
\end{equation}
If the Keller-Osserman type condition \eqref{KOgen} holds, then system \eqref{systemrad1} does not admit a radially symmetric positive entire solution. If $f,g$ also satisfy (A3), then system \eqref{systemrad1} does not admit any positive entire solutions if \eqref{KOgen} holds.
\end{cor}
In Corollary \ref{corollary} the weight functions $\alpha(|x|),\beta(|x|)$ are not assumed to satisfy conditions \eqref{inf p}-\eqref{p over q bound}, so in particular $\alpha(|x|)$ is allowed to be unbounded. Noting that $\tilde{p}(|x|)\equiv C_1$ and $\tilde{q}(|x|)\equiv C_2$ satisfy conditions \eqref{inf p}-\eqref{p over q bound} for any pair of positive constants $C_1,C_2$, if  $\alpha(|x|),\beta(|x|)$ merely satisfy $\alpha(|x|)\geq C_1>0$ and $\beta(|x|)\geq C_2>0$, then Corollary \ref{corollary} applies.  While the Keller-Osserman type condition \eqref{KOgen} is quite technical, for a wide class of nonlinearities it takes a simpler form, and is in fact sharp. Indeed, by imposing a polynomial growth condition on $g$, we have the following sharp result.

\begin{theorem}\label{main thm 3}
With the same assumptions as Theorem \ref{main thm 2}, further assume that there exists a constant $k>0$ such that $g(t)/t^{k}$ is non-increasing in $(0,\infty)$ and 
\begin{equation}\label{poly growth}
    \lim\limits_{t\to\infty}\frac{g(t)}{t^{k}}>0.
\end{equation} 
Then system \eqref{systemrad} admits radially symmetric positive entire solutions if and only if \eqref{KOcond} holds. Consequently, if (A3) is satisfied then system \eqref{systemrad} admits positive entire solutions if and only if \eqref{KOcond} holds. 
\end{theorem}

We note in particular that Theorem \ref{main thm 3} applies to system \eqref{DLSsystem}, and thus system \eqref{DLSsystem} does not admit any positive entire solutions. We emphasise again that this was previously only known in the case of radial solutions. Clearly the conditions of Theorem \ref{main thm 3} are satisfied by $g(t)=t^k, k\geq 1$, but \eqref{poly growth} allows a far more general class of nonlinearities to be considered. For instance, $g$ could be any sum of nonnegative powers of $t$ with nonnegative coefficients and highest power $t^k$, once this sum is convex. We illustrate this with the following example.

\begin{example}
\normalfont{
Consider the system
  \begin{equation}\label{ex1}
\left\{
\begin{aligned}
&\Delta u = v^k   \quad&&\mbox{ in } \bR^n,\\
&\Delta v = |\nabla u|^m \quad&&\mbox{ in } \bR^n,\\
&u>0,\ v>0 \quad&&\mbox{ in } \bR^n,
\end{aligned}
\right.
\end{equation}  
where $k,m\geq 1$. In this case, the Keller-Osserman condition \eqref{KOcond} reduces to $1\leq m\leq \frac1k$, and we conclude that system \eqref{ex1} admits solutions if and only if $k=m=1$. In particular, system \eqref{DLSsystem} where $k=1$ and $m=2$, does not admit any entire solutions. More generally, suppose we have a system of the form
\begin{equation}\label{ex2}
\left\{ \begin{aligned}
\Delta u &\geq \tilde{p}(|x|)\sum_{i=1}^{s} a_iv^{k_i}
&\quad\mbox{ in } \bR^n,\\
\Delta v &\geq \tilde{q}(|x|) f(|\nabla u|)
&\quad\mbox{ in } \bR^n,\\
u>&0,\ v>0 &\quad\mbox{ in } \bR^n,
\end{aligned} \right.
\end{equation}
  where $a_i,k_i> 0$ and $s> 1$. We further suppose that conditions (A1)-(A3) and \eqref{inf p}-\eqref{p over q bound} are satisfied. Assuming that $1\leq k_1< k_2< \cdots< k_s$, Theorem \ref{main thm 3} states that system \eqref{ex2} admits positive entire solutions if and only if the Keller-Osserman type condition \eqref{KOcond} holds for $k=k_s$. 
    }
\end{example}

\begin{remark}
\normalfont{
While we have focused on system \eqref{system} in this paper, we remark that the strategy used here is applicable to systems without gradient terms present such as
\begin{equation}\label{systemvar2}
\left\{
\begin{aligned}
&\Delta u \geq p(x)g(v)   \quad&&\mbox{ in } \bR^n, \\
&\Delta v \geq q(x)f(u) \quad&&\mbox{ in } \bR^n,\\
&u>0, v>0 \quad&&\mbox{ in } \bR^n.
\end{aligned}
\right.
\end{equation}
Theorem \ref{main thm} applies to system \eqref{systemvar2} with no modifications. When it comes to the study of semilinear systems such as \eqref{systemvar2}, there is an incredibly extensive literature. With no hope of being exhaustive, we mention the papers \cite{BVG2010,CR2002,Dam13,GS22, L2011,LW2000,T2003} and the references therein. Once again, Theorem \ref{main thm} can be sharpened for wide classes of nonlinearities using known results regarding the existence of radial solutions of the system
\begin{equation*}\label{systemvar 2}
\left\{
\begin{aligned}
&\Delta u = \tilde{p}(|x|)g(v)   \quad&&\mbox{ in } \bR^n, \\
&\Delta v = \tilde{q}(|x|)f(u) \quad&&\mbox{ in } \bR^n,\\
&u>0, v>0 \quad&&\mbox{ in } \bR^n.
\end{aligned}
\right.
\end{equation*}
}
\end{remark}

The remaining sections of this paper are organised as follows. In Section 2 we gather some preliminary results about the radial solutions of system \eqref{systemrad}, and prove our comparison result. In Section 3 we use these to prove Theorems \ref{main thm}, \ref{main thm 2} and \ref{main thm 3}.

\section{Preliminary Results and Comparison Principle}
We begin by noting that system \eqref{systemrad} admits local radially symmetric solutions in some open ball $B_R\subset\bR^n$. Indeed, assuming (A1) and (A2), a standard application of the Schauder fixed point theorem shows the problem
\begin{equation}\label{radialsys}
\begin{aligned}
    &(r^{n-1}u')'=r^{n-1}\Tilde{p}(r)g(v),\\
    &(r^{n-1}v')'=r^{n-1}\Tilde{q}(r)f(|u'|),\\
    &u>0,\quad v>0,\quad u'(0)=v'(0)=0,
\end{aligned}
\end{equation}
has a solution on some interval $[0,r_0], r_0>0$. Moreover, this holds for any pair of central values $(u(0),v(0))=(a,b)\in\bR^{+}\times \bR^{+}$. Note that system \eqref{radialsys} is the polar or radial form of \eqref{systemrad}. Suppose now that $(u,v)$ is a non-constant solution of \eqref{radialsys} on the interval $[0,R)$. By assumption, the maps $r\mapsto r^{n-1}u'$ and $r\mapsto r^{n-1}v'$ are strictly increasing for $0<r<R$, and vanish at $r=0$. It thus easily follows that $u'(r),v'(r)>0$ for all $r\in(0,R)$.
Since any local solution $(u,v)$ of system \eqref{radialsys} is increasing in both components, standard arguments show that if $(u,v)$ is not also a solution of \eqref{radialsys} on $[0,\infty)$, then $|(u(r),v(r))|\to\infty$ as $r\to R_0^{-}$, for some $0<R_0<\infty$. We now proceed with a comparison principle for sub and super radial solutions of \eqref{systemrad}, which will be the key ingredient in the proof of Theorem \ref{main thm} (cf. similar results in \cite{BGW22, DK24}). To this end, we define a \textit{super-solution} of \eqref{systemrad} in $\Omega\subset\bR^n$ to be a pair $(u,v)$ satisfying
\begin{equation*}
\left\{
\begin{aligned}
&\Delta u \geq \Tilde{p}(|x|)g(v)   \quad&&\mbox{ in } \Omega, \\
&\Delta v \geq \Tilde{q}(|x|)f(|\nabla u|) \quad&&\mbox{ in } \Omega.
\end{aligned}
\right.
\end{equation*}
\textit{Sub-solutions} are defined similarly, with both of the above inequalities reversed. By assumption (A2), any solution of \eqref{system} is automatically a super-solution of \eqref{systemrad} in $\Omega=\bR^n$.

\begin{lemma} \label{ComLemma}
Assume (A1) and (A2). Suppose that $(u,v)$ is a non-constant radial super-solution of \eqref{systemrad} in $\Omega=B_R$, $(u_0,v_0)$ is a non-constant radial sub-solution of \eqref{systemrad} in $\Omega=B_{R_0}$, and let $\overline{R}=\mbox{min}\{R,R_0\}$. Suppose also that $(u,v)$ and $(u_0,v_0)$ are positive for all
$0<r<\overline{R}$. If $u(0) > u_0(0)$ and $v(0) > v_0(0)$, then
\begin{equation} \label{comp1}
u(r)>u_0(r), \quad v(r)>v_0(r), \quad u'(r)>u_0'(r), \quad v'(r)\geq v_0'(r)
\end{equation}
throughout the whole interval $(0,\overline{R})$.
\end{lemma}

\begin{proof}
Let us denote by $(0,\Tilde{R}) \subset (0,\overline{R})$ the maximal subinterval on which
\begin{equation} \label{compp2}
u(r) > u_0(r), \qquad v(r) > v_0(r) \qquad\text{for all $0<r<\Tilde{R}$.}
\end{equation}
By assumption, the function $u$ satisfies 
\begin{equation*}
\left[ r^{n-1} u'(r) \right]' \geq r^{n-1} \Tilde{p}(r)g(v(r)),
\end{equation*}
and the function $u_0$ similarly satisfies
\begin{equation*}
\left[ r^{n-1} u_0'(r) \right]' \leq r^{n-1}\Tilde{p}(r)g(v_0(r)).
\end{equation*}
Combining these inequalities we obtain
\begin{equation*}
\left[ r^{n-1} (u(r) - u_0(r))' \right]' \geq r^{n-1}\Tilde{p}(r)\left[g(v) - g(v_0)\right].
\end{equation*}
Along the interval $(0,\Tilde{R})$, we have $v>v_0>0$ by \eqref{compp2}, so $g(v) > g(v_0)$
by (A1). The right hand side above is therefore positive, so it easily follows that $u'(r)>u_0'(r)$ for all $0<r<\Tilde{R}$.

Next, we turn to the functions $v,v_0$. Similar to above, we have
\begin{equation}\label{compp3}
\left[ r^{n-1} (v(r) - v_0(r))' \right]' \geq r^{n-1}\Tilde{q}(r)\left[f(u') - f(u'_0)\right].
\end{equation}
When it comes to the interval $(0,\Tilde{R})$, we have $u'>u_0'$ by above, so $f(u') \geq f(u_0')$, meaning the right hand side of \eqref{compp3} is non-negative. If the right hand side of \eqref{compp3} happens to be identically zero, since $v'(0)=v_0'(0)=0$, we in fact have $v'(r)=v_0'(r)$ on $(0,\Tilde{R})$. If the right hand side of \eqref{compp3} becomes positive on $(0,\Tilde{R})$, we may argue as above to conclude that $v'(r)\geq v_0'(r)$ for all $0<r<\Tilde{R}$. In both cases we see that $v-v_0$ is non-decreasing on $(0,\Tilde{R})$. 
This shows that both $u(r)-u_0(r)$ and $v(r)-v_0(r)$ are non-decreasing throughout $(0,\Tilde{R})$, and as these
functions are positive at $r=0$ by assumption, the conclusion follows.
\end{proof}
\section{Proofs of Main Results}
\subsection*{Proof of Theorem \ref{main thm}}
Suppose that $(u,v)$ is a positive entire solution of \eqref{system}. For each $r>0$, let $\overline{u}(r)$ denote the spherical average of $u$ on the ball of radius $r$, that is
\begin{equation*}
    \overline{u}(r)=\frac{1}{c_nr^{n-1}}\int_{\partial B_r}u\ d\sigma.
\end{equation*}
A basic calculation gives
\begin{equation*}
    \overline{u}'(r)=\frac{1}{c_nr^{n-1}}\int_{\partial B_r}\nabla u\cdot\frac{x}{r}\ d\sigma.
\end{equation*}
Since $\Delta u\geq p(x)g(v)>0$ for all $x\neq 0$, the divergence theorem tells us that $\overline{u}'>0$ for all $r>0$. Moreover, by the Cauchy-Schwarz inequality we have
\begin{equation}\label{avggrad}
  \overline{u}'(r)\leq \frac{1}{c_nr^{n-1}}\int_{\partial B_r}|\nabla u|\ d\sigma=\overline{|\nabla u|}.
\end{equation}
We begin by again noting that due to assumption (A2), since $(u,v)$ is a solution of \eqref{system}, it is also a super-solution of system \eqref{systemrad}, that is 
\begin{equation*}
\left\{
\begin{aligned}
&\Delta u \geq \Tilde{p}(|x|)g(v), \\
&\Delta v \geq \Tilde{q}(|x|)f\left( |\nabla u|\right).\\
\end{aligned}
\right.
\end{equation*}
Now, recalling assumptions (A1)-(A3), by taking spherical averages in the above system of inequalities, and using Jensen's inequality together with \eqref{avggrad} we obtain
\begin{equation*}
\left\{
\begin{aligned}
&\Delta \overline{u} \geq \Tilde{p}(r)g(\overline{v}), \\
&\Delta \overline{v} \geq \Tilde{q}(r)f\left( \overline{u}'\right),\\
\end{aligned}
\right.
\end{equation*}
for all $r>0$. We also have $\lim\limits_{r\to 0^{+}}\overline{u}'(r)=0$, $\lim\limits_{r\to 0^{+}}\overline{u}(r)=u(0)$ and $\lim\limits_{r\to 0^{+}}\overline{v}(r)=v(0)$. In other words, $(\overline{u},\overline{v})$ is a radially symmetric super-solution of system \eqref{systemrad}, so we are in a position to apply our comparison result Lemma \ref{ComLemma}.

Let $(a,b)\in\bR^{+}\times\bR^{+}$ satisfy $a<u(0)$ and $b<v(0)$. By the Schauder fixed point theorem, there exists a solution $(\Tilde{u},\Tilde{v})$ of \eqref{radialsys} satisfying $\Tilde{u}(0)=a$ and $\Tilde{v}(0)=b$. Applying Lemma \ref{ComLemma}, we have that $\Tilde{u}<\overline{u}$ and $\Tilde{v}<\overline{v}$ for all $r\in[0,R_0)$, where $[0,R_0)$ is the maximal interval on which $(\Tilde{u},\Tilde{v})$ solves system \eqref{radialsys}. Since $(u,v)$ is an entire solution of \eqref{system}, we must have $R_0=\infty$, which finishes the proof. $\hfill \blacksquare$

\subsection*{Proof of Theorem \ref{main thm 2}}
We begin with the following lemma. Similar results in the case that $\tilde{p},\tilde{q}$ are non-decreasing can be found in \cite{DK24, S15}. The proof here is only a slight modification, but we include it for completeness.
\begin{lemma}\label{convex}
Suppose $(u,v)$ is a non-constant positive radial solution of \eqref{systemrad} in $\Omega=B_R$. Then 
\begin{align}
     \Phi_{\tilde{p}}(r)g(v) \leq&u''(r)\leq \tilde{p}(r)g(v) ,\label{convex bounds1}\\
          \Phi_{\tilde{q}}(r)f(u')\leq&v''(r)\leq \tilde{q}(r)f(u'),\label{convex bounds2}
\end{align}
for all $0<r<R$, with $\Phi_{\tilde{p}},\Phi_{\tilde{q}}$ defined as in \eqref{phi}. In particular,  if $\tilde{p},\tilde{q}$ satisfy \eqref{inf p} and \eqref{inf q}, both $u(r)$ and $v(r)$ are convex for all $0<r<R$.
\end{lemma}
\begin{proof}
Since $u',v'>0$ for all $0<r<R$, the upper bounds in \eqref{convex bounds1} and \eqref{convex bounds2} clearly hold. When it comes to the lower bounds, we first rewrite \eqref{systemrad} as
 \begin{equation}\label{sysrad2}
	\left\{
	\begin{aligned}
	&[r^{n-1}u'(r)]'=r^{n-1}\Tilde{p}(r)g(v) &&\quad\mbox{for all}\quad 0<r<R,\\
    &[r^{n-1}v'(r)]'=r^{n-1}\Tilde{q}(r)f(u')&&\quad\mbox{for all}\quad 0<r<R,\\
    &u'(0)=v'(0)=0, u(0)> 0, v(0)>0.
	\end{aligned}
	\right.
\end{equation}
Now, since $v(r)$ is increasing for $0<r<R$, the first equation in \eqref{sysrad2} tells us
\begin{equation*}
    r^{n-1}u'(r)=\int_{0}^{r}t^{n-1}\Tilde{p}(t)g(v(t))dt
    \leq g(v(r))\int_{0}^{r}t^{n-1}\Tilde{p}(t)dt 
\end{equation*}
for all $0<r<R$. So we have the estimate 
\begin{equation}\label{u' estimate}
    \frac{n-1}{r}u'(r)\leq \frac{(n-1)g(v(r))}{r^n}\int_{0}^{r}t^{n-1}\Tilde{p}(t)dt
\end{equation}
for all $0<r<R$, which when combined with the first equation in \eqref{systemrad} yields
\begin{align*}
    u''(r)&=\Tilde{p}(r)g(v(r))-\frac{n-1}{r}u'(r)\\
    &\geq g(v(r))\left(\tilde{p}(r)-\frac{n-1}{r^n}\int_{0}^{r}t^{n-1}\Tilde{p}(t)dt\right)\\
    &= \Phi_{\tilde{p}}(r)g(v(r))
\end{align*}
whence the lower bound in \eqref{convex bounds1} holds, and $u(r)$ is convex when \eqref{inf p} holds.
Now, since we have just shown that $u'(r)$ is increasing for $0<r<R$, we can take the second equation in \eqref{sysrad2} and argue in a similar way to conclude the lower bound in \eqref{convex bounds2} also holds, and $v(r)$ is also convex when \eqref{inf q} holds. 
\end{proof}

Recall that we have defined for $s\geq 0$
\begin{equation*}
    \mathfrak{G}(s)=\int_{0}^{s}g^2(t)\ dt.
\end{equation*}
In what follows, $C$ will be used to denote a positive constant which may vary on each occurrence. In general $C=C(n,u,v,\tilde{p},\tilde{q})$, but its exact value will not be important. For simplicity we shall let $w=u'$, and we note that $w$ is positive and increasing on $(0,R)$ by Lemma \ref{convex}. 
We first note that system \eqref{systemrad} has a positive radial solution $(u,v)$ defined on some maximal interval $[0,R_{0})$. Our goal is to show that this interval is finite provided \eqref{KOgen} holds. Since it might be the case that $\tilde{q}(0)=0$, we fix $\rho\in(0,R_0)$, and recall  from \eqref{inf q}, \eqref{convex bounds1} and \eqref{convex bounds2} that we have 
    \begin{align}
        w'(r)&\leq \tilde{p}(r)g(v(r))\label{vk1}\\
        C_{\rho}f(w(r))&\leq  v''(r)\label{vk2}
    \end{align}
 for all $ r\in[\rho,R_0)$. 
 Multiplying \eqref{vk1} and \eqref{vk2} gives
    \begin{equation*}
        f(w(r))w'(r)\leq Cg(v(r))v''(r)
    \end{equation*}
    for all $ r\in[\rho,R_0)$, recalling from \eqref{p over q bound} that $\tilde{p}(r)$ is bounded. Integrating the above inequality over $[\rho,r]$ for any $r\in(\rho,R_0)$ we obtain
    \begin{align*}
        \int_{w(\rho)}^{w(r)}f(t)\ dt&\leq C\int_{\rho}^{r}g(v(t))v''(t)\ dt\\
        &\leq C\int_{0}^{r} g(v(t))v''(t)\ dt\\
        &\leq Cg(v(r))\int_{0}^{r} v''(t)\ dt\\
        &\leq Cg(v(r))v'(r),
    \end{align*}
recalling that $v'(0)=0$. Taking the extreme left and right above, one has
\begin{align}\label{wvv'}
    \int_{0}^{w(r)}f(t)dt&\leq \int_{0}^{w(\rho)}f(t)dt +Cg(v(r))v'(r)\nonumber\\ 
    &=\left(\frac{\int_{0}^{w(\rho)}f(t)dt}{g(v(r))v'(r)}+C\right)g(v(r))v'(r)\nonumber\\
    &\leq \left(\frac{\int_{0}^{w(\rho)}f(t)dt}{g(v(\rho))v'(\rho)}+C\right)g(v(r))v'(r)\nonumber\\
    &\leq  Cg(v(r))v'(r)
\end{align}
for all $r\in[\rho,R_0)$, recalling that $v'(r)$ is positive for all $r>0$. Multiplying inequalities \eqref{vk1} and \eqref{wvv'} thus gives
\begin{equation*}
    w'(r)\int_{0}^{w(r)}f(t)dt\leq Cg^2(v(r))v'(r)
\end{equation*}
for all $ r\in[\rho,R_0)$. A further integration over $[\rho,r]$ yields
\begin{equation*}
    \int_{w(\rho)}^{w(r)}\int_{0}^{s}f(t)\ dtds \leq C\int_{v(\rho)}^{v(r)}g^2(t)\ dt\leq C\int_{0}^{v(r)}g^2(t)\ dt,
\end{equation*}
for all $r\in(\rho,R_0)$. Similar to above, there exists $C>0$ such that
\begin{equation*}
   \int_{0}^{w(r)}\int_{0}^{s}f(t)\ dtds \leq C\int_{0}^{v(r)}g^2(t)\ dt
\end{equation*}
for all $r\in(\rho,R_0)$. Letting $F(s)=\int_{0}^{s}f(t)dt$ and $\mathcal{F}(s)=\int_{0}^{s}F(t)\ dt$, we may express
the above inequality as
\begin{equation*}
    \mathcal{F}(w(r)) \leq C\mathfrak{G}(v(r))
\end{equation*}
for all $r\in(\rho,R_0)$, which from \eqref{convex bounds1} yields
\begin{equation}
   \Phi_{\tilde{p}}(r) g\left(\mathfrak{G}^{-1}\left(C\mathcal{F}(w(r)) \right)\right)\leq w'(r),
\end{equation}
whence we have
\begin{equation}\label{clessw}
    \Phi_{\tilde{p}}(r)\leq \frac{w'(r)}{  g\left(\mathfrak{G}^{-1}\left(C\mathcal{F}(w(r)) \right)\right)}
.
\end{equation}
Integrating \eqref{clessw} over $[\rho,r]$ and letting $r\to R_0^{-}$ we find
\begin{align}\label{intwbounded}
    \int_{\rho}^{R_0}\Phi_{\tilde{p}}(s)\ ds&\leq \int_{w(\rho)}^{w(R_0^{-})}\frac{ds}{g\left(\mathfrak{G}^{-1}\left(C\mathcal{F}(s) \right)\right)} \nonumber\\
    &\leq \int_{w(\rho)}^{\infty}\frac{ds}{g\left(\mathfrak{G}^{-1}\left(C\mathcal{F}(s) \right)\right)} 
\end{align}
 where $w(R_0^{-})\in \bR^{+}\cup\{\infty\}$ is understood to mean $\lim_{r\to R_0^{-}}w(r)$. From \eqref{inf p}, we thus have that $R_0<\infty$ if \eqref{KOgen} holds. This proves the first part of the theorem, while the second part now follows immediately from Theorem \ref{main thm}.  $\hfill \blacksquare$
\bigskip

 \noindent To see that Corollary \ref{corollary} now follows, one need only apply the comparison principle Lemma \ref{ComLemma}. Indeed,  given a local radially symmetric solution $(u,v)$ of system \eqref{systemrad1} with central values $(u(0),v(0))$, take a local radially symmetric solution $(u_0,v_0)$ of system \eqref{systemrad} with central values satisfying $u(0)>u_0(0)$ and $v(0)>v_0(0)$. Since one has $\alpha\geq \tilde{p}$, $\beta\geq \tilde{q}$, $\phi\geq g$ and $\psi\geq f$, Lemma \ref{ComLemma} can be applied and we conclude $u(r)>u_0(r)$ and $v(r)>v_0(r)$ at every point that these functions solve their respective systems. But if the Keller-Osserman type condition \eqref{KOgen} holds, from Theorem \ref{main thm 2} one has $|(u_0,v_0)|\to\infty$ as $r\to R_0^{-}$ for some $R_0<\infty$, meaning that $|(u,v)|\to\infty$ as $r\to R^{-}$ for some $R\leq R_0$.  $\hfill \blacksquare$

\subsection*{Proof of Theorem \ref{main thm 3}} 
The first part of Theorem \ref{main thm 3} is a special case of \cite[Theorem 1.2]{D2025}, while the second part follows from Theorem \ref{main thm}. $\hfill \blacksquare$

\begin{remark}
\normalfont{
While Theorem \ref{main thm 2} tells us that condition \eqref{KOgen} is sufficient to rule out the existence of entire radial solutions of \eqref{systemrad}, an obvious question is whether or not it is also a necessary condition. Introducing the function 
\begin{equation*}
    \mathcal{G}(s)=\int_{0}^{s}\frac{g(t)}{t^{\frac12}}\ dt,
\end{equation*}
one can show that a necessary condition for the nonexistence of entire radial solutions of \eqref{systemrad} is given by
\begin{equation}\label{int sqrt f}
   \int_{}^{\infty}\frac{ds}{g\left(\mathcal{G}^{-1}\left(C\int_{0}^{s}f(t)^{\frac{1}{2}}dt\right)\right)}<\infty,
\end{equation}
for some constant $C>0$. Theorem \ref{main thm 3} tells us that condition \eqref{int sqrt f} is equivalent to condition \eqref{KOgen} when $g$ grows polynomially, but the equivalence of these conditions without a growth restriction on $g$ is an open problem.
}
\end{remark}

\section*{Acknowledgements}
The author thanks the anonymous referees for their careful reading of this manuscript, and for their suggestions which have helped to improve it. 

\section*{Funding Information}
The research conducted in this publication was funded by Research Ireland under grant number 21/FFP-A/9150.

\bibliography{bib}
\bibliographystyle{siam}
\end{document}